\documentclass[preprint,12pt]{elsarticle_mod}

\newtheorem{theorem}{Theorem}[section]

\newtheorem{lemma}[theorem]{Lemma}
\newtheorem{proposition}[theorem]{Proposition}

\newtheorem{definition}[theorem]{Definition}

\newcommand{\R}{{\mathsf{I\!R}}}

\newcommand{\Q}{{\kern.24em\vrule width.04em height1.4ex%
                 depth-.05ex\kern-.26em\mathsf Q}}
\newcommand{\C}{{\kern.24em\vrule width.04em height1.4ex%
                 depth-.05ex\kern-.26em\mathsf C}}

\usepackage{amssymb}
\journal{}

\begin{document}

\begin{frontmatter}



\title{Large-time asymptotics  of moving-reaction interfaces  involving  nonlinear Henry's law and time-dependent Dirichlet data}

\author[lab1]{Toyohiko Aiki}
\author[lab2]{Adrian Muntean}

\address[lab1]{Department of Mathematics, Faculty of Science, Tokyo Women's University, 2-8-1 Mejirodai, Bunkyo-ku, Tokyo, 112-8681, Japan, email: aikit@fc.jwu.ac.jp}
\address[lab2]{CASA - Centre for Analysis, Scientific computing and Applications, Institute for Complex Molecular Systems (ICMS), Eindhoven University of Technology, 
PO Box 513, 5600 MB,  Eindhoven, The Netherlands, email: a.muntean@tue.nl}

\begin{abstract}
We study the large-time behavior of the free boundary position capturing the one-dimensional motion of the carbonation reaction front in concrete-based materials. We extend here our rigorous justification of the $\sqrt{t}$-behavior of reaction penetration depths by including non-linear effects due to deviations from the classical Henry's law and time-dependent Dirichlet data. 

\end{abstract}


\end{frontmatter}

% \linenumbers


\section{Introduction}

In this paper, we deal with the following initial free-boundary value problem: 
Find $\{s, u, v\}$ such that 
\begin{eqnarray}
& & Q_s(T) = \{(t,x)| 0 < x < s(t), 0 < t < T\}, \nonumber \\
& &  u_t - (\kappa_1 u_x)_x = f(u,v)  \quad \mbox{ in } Q_s(T), \label{a20} \\
& &   v_t - (\kappa_2 v_x)_x = - f(u,v) \quad \mbox{ in } Q_s(T),  \label{a30} \\
& & u(t,0)= g(t), v(t,0) = h(t) \quad \mbox{ for } 0<t<T, \label{a40} \\
& & u(0,x)= u_0(x), v(0,x) = v_0(x) \quad \mbox{ for } 0 < x < s_0, \label{a50} \\ 
& & s'(t) ( = \frac{d}{dt} s(t) ) = \psi(u(t,s(t))) \quad \mbox{ for } 0<t<T,  \label{a60}\\ 
& & \kappa_1u_x(t,s(t)) = - \psi(u(t,s(t))) - s'(t) u(t,s(t)) \quad \mbox{ for } 0<t<T, 
        \label{a70}\\
& & \kappa_2 v_x(t,s(t)) = - s'(t) v(t,s(t)) \quad \mbox{ for } 0<t<T,  \label{a80} \\
& & s(0) = s_0, \label{a90} 
\end{eqnarray}
where $T > 0$, 
$\kappa_1$ and $\kappa_2$ are positive constants, $f$ is a given continuous function on 
$\mathcal \R^2$, $g$ and $h$ are boundary data, $u_0$,  $v_0$ and $s_0$ are initial data and 
 $\psi(r) =  \kappa_0|[r]^+|^p$ where $\kappa_0 > 0$ and $p \geq 1$ are given constants.    Here $u$ and $v$ represent the mass concentration of carbon dioxide dissolved in water and in air, respectively, while $s(t)$ denotes the position of the penetration reaction front in concrete at time $t>0$. The interface $s$ separates the carbonated from the non-carbonated regions.
 
   We denote by P$(f)$ the above system (\ref{a20}) $\sim$ (\ref{a90}).  P$(f)$ describes to so called {\em concrete carbonation process}, one of the most important physico-chemical mechanisms responsible for the durability  of concrete structures; see \cite{Has,Ruan} for more details of the civil engineering problem.

The target here is to study the {\em large-time behavior of weak solutions}\footnote{This is the way we translate the concept of "material durability" in mathematical terms.} in the presence of macroscopic nonlinear Henry's law and time-dependent Dirichlet boundary conditions. To get a bit the flavor of mathematical investigations of the effects by Henry's law for this or closely related reaction-diffusion systems, we refer the reader to  \cite{Ai-Mun1} (linear Henry's law) and \cite{Malte,Maria} (micro- and micro-macro Henry-like laws). 
Essentially, we are able to present refined estimates that extend the proof of the rigorous large-time asymptotics beyond the settings that we have elucidated in \cite{IFB,CPAA}. In practical terms, we show that there exist two positive constants $c_*$ and $C_*$, depending on all material parameters and initial and boundary data, such that 
\begin{equation}\label{fb}
 c_* \sqrt{t} \leq s(t) \leq C_* \sqrt{t + 1} \quad \mbox{ for } t \geq 0. 
 \end{equation}
Based on (\ref{fb}), we can now explain that the deviations of carbonation fronts from the $\sqrt{t}$-law emphasized in \cite{Kris} are certainly not due to eventual nonlinearities arising in the productions by Henry's law nor due  to the time-changing (local) atmospheric dioxide concentrations.  Therefore, there must be other reasons for this to happen. However,  we prefer  to not give rise here to many discussions in this direction. We just want to mention a first plausible reason:  Depending on the cement chemistry, the carbonation reaction might not be sufficiently fast to justify a free-boundary formulation. This fact may naturally lead to a  variety of different large-time asymptotics. 

The reminder of the paper focuses on justifying rigorously the upper and lower bounds on the interface position $s(t)$ as indicated in (\ref{fb}).

\section{Technical preliminaries. Statement of the main theorem}

We consider P$(f)$ in the cylindrical domain $Q(T):= (0,T) \times (0,1)$ by using change of variables in order to define a solution with usual notations: Let 
\begin{equation} \label{trans}
\bar{u}(t,y) = u(t,s(t)y) \mbox{ and } \bar{v}(t,y) = v(t,s(t)y)
\mbox{  for } (t,y) \in Q(T). 
\end{equation}
 Then, it holds that
  \begin{eqnarray*}
& & {\bar u}_t - \frac{\kappa_1}{s^2} \bar{u}_{yy} - \frac{s'}{s} y \bar{u}_y = f(\bar{u}, \bar{v})
 \quad \mbox{ in } Q(T), \label{eqn.a90} \\
& & {\bar v}_t - \frac{\kappa_2}{s^2} \bar{v}_{yy} - \frac{s'}{s} y \bar{v}_y = 
- f(\bar{u}, \bar{v})
 \quad \mbox{ in } Q(T), \\
& & \bar{u}(t,0) = g(t), \bar{v}(t,0) = h(t) \quad \mbox{ for } 0 < t < T, \\
& & s'(t) =  \psi(\bar{u}(t,1))  \quad \mbox{ for } 0 < t < T, \\
& & - \frac{\kappa_1}{s(t)} {\bar u}_y(t,1) = s'(t)  \bar{u}(t,1) +
  s'(t)  \quad \mbox{ for } 0 < t < T, \\
& & - \frac{\kappa_2}{s(t)} \bar{v}_y(t,1) = s'(t) \bar{v}(t, 1)
         \quad \mbox{ for } 0 < t < T, \\
& & s(0) = s_0, \bar{u}(0,y) = \bar{u}_0(y), \bar{v}(0,y) = \bar{v}_0(y) \quad \mbox{ for }
0 < y < 1, \label{eqn.a160}
\end{eqnarray*}
where $\bar{u}_0(y) = u_0(s_0 y)$ and $\bar{v}_0(y) = v_0(s_0 y)$ for $y \in [0,1]$.

For simplicity, throughout this paper we introduce the following notations related to some function spaces: 
We put $H := L^2(0,1)$, $X := \{z \in H^1(0,1): z(0) = 0\}$, $|z|_X = |z_x|_H$ for $z \in X$,
$V(T) = L^{\infty}(0, T; H) \cap L^2(0,T; H^1(0,1))$,
$V_0(T) = V(T) \cap L^2(0,T; X)$  and
$|z|_{V(T)} = |z|_{L^{\infty}(0,T; H)} + |z|_{L^2(0,T; X)}$ for $z \in V(T)$. 
Also, we denote by $X^*$ and $\langle \cdot, \cdot \rangle_X$
the dual space of $X$ and the duality pairing between $X$ and $X^*$,
respectively.

By using these notations we define a weak solution of P$(f)$  in the following way: 

\begin{definition} \label{weak_solution}
Let $s$ be a function on $[0,T]$ and $u$, $v$ be functions on $Q_s(T)$
for $0 < T < \infty$. 
We call that a triplet $\{s, u, v\}$ is a weak solution of P$(f)$ on $[0,T]$ 
if the conditions (S1) $\sim$ (S5) hold: \\
(S1) $s \in W^{1,\infty}(0,T)$ with $s > 0$ on $[0,T]$, $(\bar{u},\bar{v}) \in (W^{1,2}(0,T; X^*) \cap V(T) \cap
L^{\infty}(Q(T)))^2$.
\\
(S2) $\bar{u} - g, \bar{v} - h \in L^2(0,T; X)$, $s(0) = s_0$, 
$u(0) = u_0$ and $v(0) =  v_0$ on $[0, s_0]$. 
\\
(S3) $s'(t) = \psi(u(t,s(t))$ for a.e. $t \in [0,T]$. 
$$
\int_0^{T} \langle \bar{u}_t(t), z(t) \rangle_X dt
+ \int_{Q(T)} \frac{\kappa_1}{s^2(t)} \bar{u}_y(t) z_y(t) dy dt
+ \int_0^{T}
 (\frac{s'(t)}{s(t)} \bar{u}(t,1) + \frac{s'(t)}{s(t)} )z(t,1) dt  \leqno{\mbox{(S4)}}
 $$
$$ = \int_{Q(T)} (f(\bar{u}(t),\bar{v}(t)) + \frac{s'(t)}{s(t)} y \bar{u}_y(t)) z(t) dydt  
 \quad \mbox{ for }
 z \in V_0(T).
$$
$$
\int_0^{T} \langle \bar{v}_t(t), z(t) \rangle_X dt
+ \int_{Q(T)} \frac{\kappa_2}{s^2(t)} \bar{v}_y(t) z_y(t) dy dt
+ \int_0^{T}
 \frac{s'(t)}{s(t)} \bar{v}(t,1)  z(t,1) dt  \leqno{\mbox{(S5)}}
 $$
$$ = \int_{Q(T)} (- f(\bar{u}(t),\bar{v}(t)) + \frac{s'(t)}{s(t)} y \bar{v}_y(t)) z(t) dy dt  
 \quad \mbox{ for }
 z \in V_0(T).
$$
Moreover, let $s$ be a function on $(0,\infty)$, and $u$ and $v$ be functions on $Q_s := \{ 
(t,x)| t > 0, 0 < x < s(t)\}$. We say that $\{s, u, v\}$ is a weak solution of P$(f)$ on $[0,\infty)$ if for any $T > 0$ the triplet $\{s, u, v\}$  is a weak solution of P$(f)$ on $[0,T]$. 
\end{definition}

Next, we give a list of assumptions for data as follows:

(A1) $f(u, v) = \phi(\gamma v - u)$ and $\phi$ is a locally Lipschitz 
continuous and increasing function on $\R$ with $\phi(0) = 0$ and  
$$ \phi(r)r  \geq C_{\phi} |r|^{1+q} \mbox{ for } r \in {\R},  $$ 
where $q \geq 1$ and $C_{\phi}$ is a positive constant. 
  
(A2) $g, h \in W^{1,2}_{loc}([0,\infty)) \cap L^{\infty}(0,\infty)$, 
$0 < g_0 \leq g$, $h \geq 0$ on $[0,\infty)$,  and 
$g - g_*, h - h_* \in W^{1,1}(0,\infty)$, where $g_0$, $g_*$ and $h_*$ are positive constants with 
$\gamma h_* = g_*$.

(A3) $s_0 > 0$ and $u_0, v_0 \in L^{\infty}(0,s_0)$,  $u_0, v_0 \geq 0$ a.e. on $(0,s_0)$.

\vskip 12pt
Our main result is as follows: 

\begin{theorem} \label{main}
If (A1), (A2) and (A3) hold, then the problem P$(f)$ has a weak solution $\{s, u, v\}$ on $[0,\infty)$. Moreover, there exist two positive constants $c_*$ and $C_*$ such that 
$$ c_* \sqrt{t} \leq s(t) \leq C_* \sqrt{t + 1} \quad \mbox{ for } t \geq 0. $$
\end{theorem}

\vskip 12pt
In order to prove Theorem \ref{main} 
we introduce the following notations: For $m > 0$ we put 
$$ \phi_m(r) = \left\{ \begin{array}{ll}
           \phi(m) & \mbox{ for } r > m, \\
           \phi(r) & \mbox{ for } |r| \leq m, \\ 
           \phi(-m) & \mbox{ for } r < -m, \end{array} \right.  \quad
$$
and $f_m(u, v) = \phi_m(\gamma v - u)$ for $(u,v) \in \R^2$. 
Obviously, for each  $m  > 0$ $\phi_m$ and $f_m$ are Lipschitz continuous. 
Then, we can denote by $C_m$ the common Lipschitz constant of $\phi_m$ and $f_m$.

Let $s \in W^{1,2}(0,T)$ and $m > 0$.  By using these notations we 
consider the auxiliary  problem SP$_m(s, \bar{u}_0, \bar{v}_0$):$= \{$(\ref{1-11}) $\sim$ (\ref{1-12})$\}$.  
\begin{eqnarray}
& & {\bar u}_t - \frac{\kappa_1}{s^2} \bar{u}_{yy} - \frac{s'}{s} y \bar{u}_y = 
 f_m(\bar{u}, \bar{v})
 \quad \mbox{ in } Q(T), \label{1-11} \\
& & {\bar v}_t - \frac{\kappa_2}{s^2} \bar{v}_{yy} - \frac{s'}{s} y \bar{v}_y = 
- f_m(\bar{u}, \bar{v})
 \quad \mbox{ in } Q(T), \\
& & \bar{u}(t,0) = g(t), \bar{v}(t,0) = h(t) \quad \mbox{ for } 0 < t < T, \\
& & - \frac{\kappa_1}{s(t)} {\bar u}_y(t,1) = s'(t)  \bar{u}(t,1) +
  \psi(\bar{u}(t,1))  \quad \mbox{ for } 0 < t < T, \\
& & - \frac{\kappa_2}{s(t)} \bar{v}_y(t,1) = s'(t) \bar{v}(t, 1)
         \quad \mbox{ for } 0 < t < T, \\
& & \bar{u}(0,y) = \bar{u}_0(y), \bar{v}(0,y) = \bar{v}_0(y) \quad \mbox{ for }
0 < y < 1, \label{1-12}
\end{eqnarray}
where $\bar{u}_0$ and $\bar{v}_0$ are given functions on the interval $[0,1]$. 

Relying on the basic properties of the  solutions to SP$_m(s, \bar{u}_0, \bar{v}_0$) (as indicated in the next 
section),  we will be able  prove our main result, that is Theorem \ref{main}, in the last section of the paper. 

\section{Basic results for the auxiliary problem SP$_m(s, \bar{u}_0, \bar{v}_0$)} \label{lemmas} 

We begin the section by showing a first result concerned with the solvability for the problem 
SP$_m(s, \bar{u}_0, \bar{v}_0$). 

\begin{proposition}\label{ex}
Let $m > 0$, $T > 0$, $s \in W^{2,1}(0, T)$ with $s(0) > 0$ and $s' \geq 0$ on $[0,T]$, 
$g, h \in W^{1,2}(0,T)$, $\bar{u}_0 - g(0) \in X$ and $\bar{v}_0 - h(0) \in X$. 
Then there exist one and only one pair $(\bar{u}, \bar{v}) \in (W^{1,2}(0,T; H) \cap L^{\infty}(0,T; H^1(0,1)) \cap L^2(0,T; H^2(0,1)))^2$ satisfying (\ref{1-11}) $\sim$ (\ref{1-12}) in the usual sense, that is, $(\bar{u}, \bar{v})$ is a unique solution of SP$_m(s, \bar{u}_0, \bar{v}_0$) on $[0,T]$. 
\end{proposition}

We can prove this proposition in a way quite similar to the working strategy illustrated in  the proofs from Section 2 in  \cite{Ai-Mun1}. Essentially, we rely on a  Banach's fixed point argument. We omit here the proof and refer the reader to \cite{Ai-Mun1}. 

As next step, we establish the positivity and the existence of upper bounds for a solution of SP$_m(s, \bar{u}_0, \bar{v}_0$). 

\begin{lemma} \label{bounds}
Under the same assumptions as in Proposition \ref{ex} let  $(\bar{u}, \bar{v})$ be a solution of 
SP$_m(s, \bar{u}_0, \bar{v}_0$) on $[0,T]$. 
If $0 \leq \bar{u}_0 \leq u_*$ and $0 \leq \bar{v}_0 \leq v_*$ on $[0,1]$,  
 $0 \leq g \leq u_*$ and $0 \leq h \leq v_*$ on $[0,T]$ and $u_* = \gamma v_*$, where $u_*$ and $v_*$ are positive constants, then 
$$ 0 \leq \bar{u} \leq u_*,  0 \leq \bar{v} \leq v_* \mbox{ on } Q(T). $$ 
\end{lemma}

{\it Proof. } 
We multiply (\ref{1-11}) by $-[-\bar{u}]^+$ to obtain

\begin{eqnarray*}
& & \frac{1}{2} \frac{d}{dt} |[-\bar{u}]^+|_H^2 
+ \frac{\kappa_1}{s^2} \int_0^1 |[-\bar{u}]_y^+|^2 dy 
- \frac{s'}{s} \bar{u}(\cdot,1) [- \bar{u}(\cdot, 1)]^+ \\
& & - \frac{1}{s} \psi(\bar{u}(\cdot,1)) [- \bar{u}(\cdot, 1)]^+  \\
& = & - \int_0^1 \phi_m( \gamma \bar{v} - \bar{u}) [- \bar{u}]^+ dy 
- \frac{s'}{s} \int_0^1 y \bar{u}_y [- \bar{u}]^+ dy 
\quad \mbox{ a.e. on } [0,T]. 
\end{eqnarray*}
Here, we note that 
\begin{eqnarray*}
 -  \phi_m( \gamma \bar{v} - \bar{u}) [- \bar{u}]^+
& \leq & -  \phi_m( - \gamma [- \bar{v}]^+ - \bar{u}) [- \bar{u}]^+  \\
& \leq & C_m( \gamma [- \bar{v}]^+ + |\bar{u}|) [- \bar{u}]^+  \\
& \leq & C_m(\gamma + 1) (|[- \bar{v}]^+ [- \bar{u}]^+ + |[- \bar{u}]^+|^2)
\quad \mbox{ a.e. on }Q(T), 
\end{eqnarray*}
and 
$$   \psi(\bar{u}(\cdot,1)) [- \bar{u}(\cdot, 1)]^+ = 0 \quad \mbox{ a.e. on }Q(T). 
$$
Then, it follows that 
$$
   \frac{1}{2} \frac{d}{dt} |[-\bar{u}]^+|_H^2 
+ \frac{\kappa_1}{2s^2} |[-\bar{u}]_y^+|_H^2    
 \leq  C_{1m}  (|[- \bar{v}]^+|_H^2 + |[- \bar{u}]^+|_H^2) 
 \mbox{ a.e. on } [0,T], 
$$
where $C_{1m} = 2C_m(\gamma + 1) +  \frac{1}{2\kappa_1} |s'|_{L^{\infty}(0,T)}^2$. 

Similarly, we can show that 
$$   \frac{1}{2} \frac{d}{dt} |[-\bar{v}]^+|_H^2 
+ \frac{\kappa_2}{2s^2} |[-\bar{v}]_y^+|_H^2    \\
\leq   C_{2m}  (|[- \bar{v}]^+|_H^2 + |[- \bar{u}]^+|_H^2)
 \mbox{ a.e. on } [0,T],  
$$
where $C_{2m} = 2C_m(\gamma + 1) +  \frac{1}{2\kappa_2} |s'|_{L^{\infty}(0,T)}^2$. From the above inequalities Gronwall's lemma implies that $[- \bar{u}]^+ =0$ and $[- \bar{v}]^+= 0$ a.e on $Q(T)$, that is, $\bar{u} \geq 0$ and $\bar{v} \geq 0$ a.e. on $Q(T)$. 

From now on we shall show the boundedness of the solutions. First, by (A2) and (A3) we can take positive constants $u_*$ and $v_*$ satisfying the inequality in the assumption of this Lemma. 

Next, we multiply (\ref{1-11}) by $[\bar{u} - u_*]^+$ and have
\begin{eqnarray*}
& & \frac{1}{2} \frac{d}{dt} |[\bar{u} - u_*]^+|_H^2 
+ \frac{\kappa_1}{s^2} |[\bar{u}- u_*]_y^+|_H^2  
+ \frac{s'}{s} \bar{u}(\cdot,1) [ \bar{u}(\cdot, 1) - u_*]^+ \\
& & + \frac{1}{s} \psi(\bar{u}(\cdot,1)) [ \bar{u}(\cdot, 1) - u_*]^+  \\
& = &   \int_0^1 \phi_m( \gamma \bar{v} - \bar{u}) [ \bar{u} - u_*]^+ dy 
+ \frac{s'}{s} \int_0^1 y \bar{u}_y [ \bar{u} - u_*]^+ dy 
\quad \mbox{ a.e. on } [0,T]. 
\end{eqnarray*}

Similarly, we see that 
\begin{eqnarray*}
& & \frac{1}{2} \frac{d}{dt} |[\bar{v} - v_*]^+|_H^2 
+ \frac{\kappa_2}{s^2} |[\bar{v}- v_*]_y^+|_H^2  
+ \frac{s'}{s} \bar{v}(\cdot,1) [ \bar{v}(\cdot, 1) - v_*]^+ \\
& = &  - \int_0^1 \phi_m( \gamma \bar{v} - \bar{v}) [ \bar{v} - v_*]^+ dy 
+ \frac{s'}{s} \int_0^1 y \bar{v}_y [ \bar{v} - v_*]^+ dy 
\quad \mbox{ a.e. on } [0,T]. 
\end{eqnarray*}
Here, elementary calculations lead to 
\begin{eqnarray*}
& & \phi_m(\gamma\hat{v} - \bar{u}) ([ \bar{u} - u_*]^+  -  [ \bar{v} - v_*]^+ ) \\
& = & \phi_m(\gamma(\hat{v} - v_*) - (\bar{u} - u_*)) ([ \bar{u} - u_*]^+  -  [ \bar{v} - v_*]^+ ) 
\\
& \leq & \phi_m(\gamma([\hat{v} - v_*]^+) - (\bar{u} - u_*)) [ \bar{u} - u_*]^+  
- \phi_m(\gamma(\hat{v} - v_*) - [\bar{u} - u_*]^+)[\bar{v} - v_*]^+\\
& \leq & C_{3m}(|[\bar{u} - u_*]^+|^2 + |[\hat{v} - v_*]^+|^2 )
 \quad \mbox{ a.e. on } Q(T),
\end{eqnarray*}
where $C_{3m} = 2C_m \gamma + C_m( \gamma + 1)$. 

From the above inequalities it follows that 
\begin{eqnarray*}
& & \frac{1}{2} \frac{d}{dt} (|[\bar{u} - u_*]^+|_H^2 +  |[\bar{v} - v_*]^+|_H^2 ) 
 + \frac{\kappa_1}{s^2}  |[\bar{u}- u_*]_y^+|_H^2 
+ \frac{\kappa_2}{s^2}  |[\bar{v}- v_*]_y^+|_H^2 \\
& \leq & 
C_{3m} (|[\bar{u} - u_*]^+|_H^2 + |[\hat{v} - v_*]^+|_H^2 ) \\
& & + \frac{\kappa_1}{2s^2}  |[\bar{u}- u_*]_y^+|_H^2 
+ \frac{1}{2\kappa_1} |s'|_{L^{\infty}(0,T)}^2 |[\bar{u}- u_*]^+|_H^2 \\
& & + \frac{\kappa_2}{2s^2}  |[\bar{v}- v_*]_y^+|_H^2 
+ \frac{1}{2\kappa_2} |s'|_{L^{\infty}(0,T)}^2 |[\bar{v}- v_*]^+|_H^2 
\quad \mbox{ a.e. on } [0,T]
\end{eqnarray*}
so that
\begin{eqnarray*}
& & \frac{1}{2} \frac{d}{dt} 
 (|[\bar{u} - u_*]^+|_H^2 +  |[\bar{v} - v_*]^+|_H^2 )  \\
& \leq & 
 (C_{3m} + \frac{1}{2} |s'|_{L^{\infty}(0,T)}^2 
 (\frac{1}{2\kappa_1} + \frac{1}{2\kappa_2}) 
            (|[\bar{u} - u_*]^+|_H^2 +  |[\bar{v} - v_*]^+|_H^2 ) 
   \quad \mbox{ a.e. on } [0,T].
\end{eqnarray*}
Hence, by applying Gronwall's lemma we conclude that $\bar{u} \leq u_*$ and $\bar{v} \leq v_*$ 
a.e. on $Q(T)$. Thus we have proved this lemma. \hfill $\Box$

\begin{lemma} \label{ineq}
Under the same assumptions as in Proposition \ref{ex} let  $(\bar{u}, \bar{v})$ be a solution of 
SP$_m(s, \bar{u}_0, \bar{v}_0$) on $[0,T]$. If $u(t,x) = \bar{u}(t, \frac{x}{s(t)})$ and 
$v(t,x) = \bar{v}(t, \frac{x}{s(t)})$ for $(t,x) \in Q_s(T)$, then the following inequality holds:
\begin{eqnarray}
& & \frac{1}{2} \frac{d}{dt} \int_0^{s} |u - g|^2 dx 
+  \frac{\gamma}{2} \frac{d}{dt} \int_0^{s} |v - h|^2 dx \nonumber \\
& & + \kappa_1 \int_0^{s} |u_x|^2 dx 
  + \kappa_2 \gamma \int_0^{s} |v_x|^2 dx 
+ \psi(u(\cdot,s)) (u(\cdot,s) - g)    \nonumber
 \\
& & + \frac{1}{2} s'( |u(\cdot,s) - g|^2 +  \gamma |v(\cdot,s) - h|^2) 
+ C_{\phi} \int_0^{s} |\gamma v - u|^{q + 1} dx \nonumber \\
& \leq &  
- g' \int_0^{s} (u - g)dx 
- \gamma h' \int_0^{s} (v - h)dx \label{2x} \\
& & -  \int_0^{s} \phi_m (\gamma v - u)  (g - g_*) dx  
+ \gamma \int_0^{s} \phi_m (\gamma v - u)  (h - h_*) dx \nonumber \\
& & - s' g (u(\cdot,s) - g) - \gamma s' h (v(\cdot,s) - h) 
\quad \mbox{  a.e. on } [0,T]. \nonumber
\end{eqnarray}
\end{lemma}
{\it Proof. } Since $(\bar{u}, \bar{v})$ is a strong solution of SP$_m(s, \bar{u}_0, \bar{v}_0$), it holds that 
\begin{eqnarray}
& & u_t - \kappa_1 u_{xx} =   f_m(u, v)
 \quad \mbox{ in } Q(T), \label{1-21} \\
& &  v_t -  \kappa_2  v_{xx}  = - f_m(u, v)
 \quad \mbox{ in } Q(T),\label{1-22} \\
& & u(0,t) = g(t), v(0,t) = h(t) \quad \mbox{ for } 0 < t < T, \nonumber \\
& & - \kappa_1 u_x(t,s(t)) = s'(t)  u(t,s(t)) +   \psi(u(t,s(t)))  
             \quad \mbox{ for } 0 < t < T, \nonumber \\
& & - \kappa_2 v_x(t,s(t)) = s'(t) v(t, s(t))  
  \quad \mbox{ for } 0 < t < T, \nonumber \\
& &  u(0,x) = u_0(x), v(0,x) = v_0(x) \quad \mbox{ for } 0 < x < s_0. \nonumber
\end{eqnarray}
Here, we multiply (\ref{1-21}) by $u - g$ 
 and  (\ref{1-22}) by $\gamma(v - h)$,  and by using integration by parts and the boundary 
conditions we obtain 
\begin{eqnarray*}
& & \frac{1}{2} \frac{d}{dt} \int_0^{s(t)} |u(t) - g(t)|^2 dx 
+ \kappa_1 \int_0^{s(t)}|u_x(t)|^2 dx \\
& & + s'(t) |u(t,s(t))- g(t)|^2
+ \psi(u(t,s(t)))  (u(t,s(t))- g(t))
\\
& = & - s'(t) g(t) (u(t,s(t)) - g(t)) - g'(t) \int_0^{s(t)} (u(t) - g(t)) dx \\
& &   + \int_0^{s(t)} f_m(u(t), v(t))  (u(t) - g(t))dx \quad 
\mbox{ for a.e. } t \in [0,T] 
\end{eqnarray*}
and 
\begin{eqnarray*}
& & \frac{\gamma}{2} \frac{d}{dt} \int_0^{s(t)} \!\!|v(t) - g(t)|^2 dx 
+ \gamma \kappa_2 \int_0^{s(t)} \!\!|v_x(t)|^2 dx  \\
&&+ \gamma s'(t) |v(t,s(t))- h(t)|^2
\\
& = &  - \gamma  s'(t) h(t) (v(t,s(t)) - h(t)) - \gamma h'(t) \int_0^{s(t)} (v(t) - h(t)) dx \\
& &   - \gamma \int_0^{s(t)} f_m(u(t), v(t))  (v(t) - h(t))dx  \quad 
\mbox{ for a.e. } t \in [0,T]. 
\end{eqnarray*}
It is easy to see that 
\begin{eqnarray*}
& & f_m(u,v) \{(u - g) - \gamma (v - h)\} \\ 
& = &  - \phi_m(\gamma v - u) (\gamma v - u) - \phi_m(\gamma v - u) \{ g - g_* + \gamma (h_* - h(t)\} \\
& \leq & - C_{\phi} |\gamma v - u|^{q+1} - \phi_m(\gamma v - u) (g - g_*)
-  \gamma \phi_m(\gamma v - u) (h_* - v) \mbox{ a.e. on } Q_s(T).  
\end{eqnarray*}
Combining these inequalities leads in a straightforward manner to the conclusion of this Lemma.  \hfill $\Box$

\vskip 12pt
The aim of this section is to establish the existence and the uniqueness of a weak solution of 
SP$_m(s, \bar{u}_0, \bar{v}_0$) in case $s \in W^{1,4}(0,T)$. 
Here, we define a weak solution of SP$_m(s, \bar{u}_0, \bar{v}_0)$

\begin{definition}\label{def3}
Let  $\bar{u}$, $\bar{v}$ be functions on $Q(T)$ for $0 < T < \infty$. 
We call that a pair $\{\bar{u}, \bar{v}\}$ is a weak solution of SP$_m(s, \bar{u}_0, \bar{v}_0$) on $[0,T]$ 
if the conditions (SS1) $\sim$ (SS4) hold: \\
(SS1)  $(\bar{u},\bar{v}) \in (W^{1,2}(0,T; X^*) \cap V(T) \cap L^{\infty}(Q(T)))^2$.
\\
(SS2) $\bar{u} - g, \bar{v} - h \in L^2(0,T; X)$, $\bar{u}(0) = \bar{u}_0$ and 
$\bar{v}(0) =  \bar{v}_0$.
\\
$$
\int_0^{T} \langle \bar{u}_t, z \rangle_X dt
+ \int_{Q(T)} \frac{\kappa_1}{s^2} \bar{u}_y z_y dy dt
+ \int_0^{T}
 (\frac{s'}{s} \bar{u}(\cdot,1) + \frac{1}{s} \psi(\bar{u}(\cdot,1))) z(\cdot,1) dt  
   \leqno{\mbox{(SS3)}}
 $$
$$ = \int_{Q(T)} (f_m(\bar{u},\bar{v}) + \frac{s'}{s} y \bar{u}_y) z dydt  
 \quad \mbox{ for }
 z \in V_0(T).
$$
$$
\int_0^{T} \langle \bar{v}_t, z \rangle_X dt
+ \int_{Q(T)} \frac{\kappa_2}{s^2} \bar{v}_y z_y dy dt
+ \int_0^{T}
 \frac{s'}{s} \bar{v}(\cdot,1)  z(\cdot,1) dt  \leqno{\mbox{(SS4)}}
 $$
$$ = \int_{Q(T)} (- f_m(\bar{u},\bar{v}) + \frac{s'}{s} y \bar{v}_y) z dy dt  
 \quad \mbox{ for }
 z \in V_0(T).
$$
\end{definition}

\begin{proposition}\label{prop2}
Let $T > 0$, $m > 0$, $s \in W^{1,4}(0,T)$ with $s(0) > 0$, $s' \geq 0$ a.e. on $[0,T]$, 
$g, h \in W^{1,2}(0,T)$ with $g, h \geq 0$ on $[0,T]$ and 
$\hat{u}_0, \hat{v}_0 \in L^{\infty}(0,1)$ with $\hat{u}_0, \hat{v}_0 \geq 0$ 
a.e. on $[0,1]$. 
Then SP$_m(s, \bar{u}_0, \bar{v}_0$) has a unique weak solution 
$\{\bar{u}, \bar{v}\}$  on $[0,T]$. 
Moreover, (\ref{2x}) holds a.e. on $[0,T]$ with $\{u, v\}$, where 
$u(t,x) = \bar{u}(t,\frac{x}{s(t)})$ and 
$v(t,x) = \bar{v}(t,\frac{x}{s(t)})$ for $(t,x) \in Q_s(T)$.

\end{proposition}

{\it Proof. } 
First, we take sequences $\{s_n\} \subset W^{2,1}(0,T)$, $\{\bar{u}_{0n}\} \subset H^1(0,1)$ and 
$\{\bar{v}_{0n}\} \subset H^1(0,1)$ such that 
$s_n \to s$ in $W^{1,4}(0,T)$ as $n \to \infty$, $s_n(0) = s(0)$, $s_n' \geq 0$ on $[0,T]$ for $n$, 
$\bar{u}_{0n} \to \bar{u}_0$ and $\bar{v}_{0n} \to \bar{v}_0$ in $H$ as $n \to \infty$, 
$0 \leq \bar{u}_{0n} \leq |\bar{u}_{0}|_{L^{\infty}(0,1)} + 1$,  
$0 \leq \bar{v}_{0n} \leq |\bar{v}_{0}|_{L^{\infty}(0,1)} + 1$ on $[0.1]$ and 
$\bar{u}_{0n} - g(0), \bar{v}_{0n} - h(0) \in X$ for $n$. 
Obviously, there exists a positive constant $L$ such that $0 < s(0) \leq s_n \leq L$ on $[0,T]$ for $n$. 
  Then, Proposition \ref{ex} and Lemma \ref{bounds} imply that SP$_m(s_n, \bar{u}_{0n}, \bar{v}_{0n}$) has a solution $(\bar{u}_n, \bar{v}_n)$ on $[0,T]$ and 
$0 \leq \bar{u}_n \leq u_*$ and $0 \leq \bar{v}_n \leq v_*$ on $Q(T)$ for each $n$, where $u_*$ and $v_*$ are positive constants satisfying 
$u_* \geq \max\{ |\bar{u}_{0}|_{L^{\infty}(0,1)} + 1, |g|_{L^{\infty}(0,T)}\}$, 
$v_* \geq \max\{ |\bar{v}_{0}|_{L^{\infty}(0,1)} + 1, |h|_{L^{\infty}(0,T)}\}$ and 
$u_* = \gamma v_*$. Moreover, by Lemma \ref{ineq} and putting 
$u_n(t,x) = \bar{u}_n(t, \frac{x}{s_n(t)})$ and 
$v_n(t,x) = \bar{v}_n(t, \frac{x}{s_n(t)})$ for $(t,x) \in Q_{s_n}(T)$,  we see that 
\begin{eqnarray*}
& & \frac{1}{2} \frac{d}{dt} \int_0^{s_n} |u_n - g|^2 dx 
+  \frac{\gamma}{2} \frac{d}{dt} \int_0^{s_n} |v_n - h|^2 dx \nonumber \\
& & + \kappa_1 \int_0^{s_n} |u_{nx}|^2 dx 
  + \kappa_2 \gamma \int_0^{s_n} |v_{nx}|^2 dx 
+ \psi(u_n(\cdot, s_n)) ( u_n(\cdot, s_n) - g)   \nonumber
 \\
& & + \frac{1}{2} s_n'( |u_n(\cdot,s_n)) - g|^2 +  \gamma |v_n(\cdot,s_n) - h|^2) 
+ C_{\phi} \int_0^{s_n} \!\!|\gamma v_n - u_n|^{q + 1} dx \nonumber \\
& \leq & 
- g' \int_0^{s_n} (u_n - g)dx 
- \gamma h' \int_0^{s_n} (v_n - h)dx \label{2q} \\
& & -  \int_0^{s_n} \phi_m (\gamma v_n - u_n)  (g - g_*) dx  
+ \gamma \int_0^{s_n} \phi_m (\gamma v_n - u_n)  (h - h_*) dx \nonumber \\
& & - s_n' g (u_n(\cdot,s_n) - g) 
- \gamma s_n' h (v_n(\cdot,s_n) - h) 
\quad \mbox{  a.e.  on } [0,T]. \nonumber
\end{eqnarray*}
Here, we note that
$$ |\phi_m(\gamma v_n - u_n)| \leq \phi(\gamma v_*) - \phi(-u_*) =: C_4 \quad \mbox{ on } Q_{s_n}(T),  $$
$$ \psi(u_n(\cdot, s_n)) ( u_n(\cdot, s_n) - g) \geq 
 \hat{\psi}(u_n(\cdot, s_n)) - \hat{\psi}(g) \quad \mbox{ a.e. on } [0,T], $$ 
where $\hat{\psi}(r) = \int_0^r \psi(\xi) d\xi$ for $ r \in \R$.  
Then by using Young's inequality we obtain 
\begin{eqnarray*}
& & \frac{1}{2} \frac{d}{dt} \int_0^{s_n} |u_n - g|^2 dx 
+  \frac{\gamma}{2} \frac{d}{dt} \int_0^{s_n} |v_n - h|^2 dx \nonumber \\
& & + \kappa_1 \int_0^{s_n} |u_{nx}|^2 dx 
  + \kappa_2 \int_0^{s_n} |v_{nx}|^2 dx 
+ \hat{\psi}(u_n(\cdot,s_n))   \nonumber
 \\
& & + \frac{1}{4} s_n'( |u_n(\cdot,s_n) - g|^2 +  \gamma |v_n(\cdot,s_n) - h|^2) 
+ C_{\phi} \int_0^{s_n} \!\!|\gamma v_n - u_n|^{q + 1} dx \nonumber \\
& \leq & \hat{\psi}(g) 
+ \frac{1}{2} (|g'|^2 + \gamma |h'|^2 )+ 
\frac{1}{2}( \int_0^{s_n} |u_n - g|^2 dx + \gamma  \int_0^{s_n} |v_n - h|^2 dx) \label{2k} \\
& & + C_4  \int_0^{s_n}  (|g - g_*| + \gamma  |h - h_*|) dx 
 + s_n' (|g|^2 + \gamma |h|^2)  
\quad \mbox{  a.e.  on } [0,T]. \nonumber
\end{eqnarray*}
Hence, by applying Gronwall's lemma we observe that 
$$  \int_0^T \int_0^{s_n} (|u_{nx}|^2   +  |v_{nx}|^2) dxdt \leq C \mbox{ for } n,
$$
where $C$ is a positive constant independent of  $n$. This implies that 
$\{\bar{u}_n\}$ and $\{\bar{v}_n\}$ are bounded in $L^2(0,T; H^1(0,1))$, since 
$|\bar{u}_{ny}(t)|^2_H  = s_n(t) \int_0^{s_n(t)} |u_{nx}(t)|^2 dx$ for $t \in [0,T]$. 

Next, we provide the boundedness of $\bar{u}_{nt}$ and $\bar{v}_{nt}$. 
Let $\eta \in X$. Then it is easy to see that 
\begin{eqnarray}
& & |\int_0^1 \bar{u}_{nt} \eta dy| \nonumber \\
& = & 
|\int_0^1 (\frac{\kappa_1}{s_n^2} \bar{u}_{nyy} + \frac{s_n'}{s_n} y \bar{u}_{ny} + 
       f_m(\bar{u}_n, \bar{v}_n)) \eta dy | \nonumber \\
& \leq &  \frac{\kappa_1}{s_n^2} \int_0^1 |\bar{u}_{ny}||\eta_y| dy 
+ \frac{s_n'}{s_n}|\bar{u}_n(\cdot,1) \eta(1)| + |\psi(\bar{u}_n(\cdot,1)) \eta(1)|
\nonumber  \\
& &  + \frac{s_n'}{s_n} |\int_0^1  \bar{u}_{n} (\eta + y \eta_y) dy | 
   + \frac{s_n'}{s_n} |\bar{u}_{n}(\cdot,1) \eta(1) | 
 + |\int_0^1 f_m(\bar{u}_n, \bar{v}_n)) \eta dy | \nonumber \\
& \leq &  \frac{\kappa_1}{s_n^2} |\bar{u}_{ny}|_H |\eta_y|_H 
+ \frac{s_n'}{s_n} u_* |\eta(1)| + |\psi(u_*)|| \eta(1)| \label{22} \\
& &  + \frac{s_n'}{s_n} u_* (|\eta|_H + |\eta_y|_H) 
   + \frac{s_n'}{s_n} u_* |\eta(1)| 
 + C_4 |\eta|_H \quad \mbox{ a.e on }  [0,T]. \nonumber
\end{eqnarray}
On account of the boundedness of $\{\bar{u}_{n}\}$ in $L^2(0,T; H^1(0,1))$ we infer that 
$\{ \bar{u}_{nt} \}$ is bounded in $L^2(0,T; X^*)$. 
Similarly,  $\{ \bar{v}_{nt} \}$ is also bounded in $L^2(0,T; X^*)$.

From the above uniform estimates there exist a subsequence $\{n_j\} \subset \{n\}$ and $(\bar{u}, \bar{v})$ such that  $(\bar{u}, \bar{v})$ satisfies (SS1), 
$\bar{u}_{n_j} \to \bar{u}$ and $\bar{v}_{n_j} \to \bar{v}$ weakly* in $L^{\infty}(Q(T))$, 
weakly in $L^2(0,T; H^1(0,1))$ and weakly in  $W^{1,2}(0,T; X^*)$ as $j \to \infty$. 
Accordingly, Aubin's compactness theorem (see \cite{Lions}) implies that  
$\bar{u}_{n_j} \to \bar{u}$ and $\bar{v}_{n_j} \to \bar{v}$ in $L^{2}(0,T;  H)$ as $j \to \infty$. 
Moreover, $\bar{u}_{n_j}(t) \to \bar{u}(t)$ and $\bar{v}_{n_j}(t) \to \bar{v}(t)$ weakly in $H$ as $j \to \infty$ for any $t \in [0,T]$, (SS2) is valid, and 
$0 \leq \bar{u} \leq u_*$ and $0 \leq \bar{v} \leq v_*$ a.e. on $Q(T)$. 

Now, I shall prove (SS3). Let $z \in V_0(T)$. Then it holds that 
\begin{eqnarray}
& & \!\!\!\!  \int_0^{T}\!\!\!\int_0^1 \bar{u}_{nt} z dx dt
+\! \int_{Q(T)} \!\!\frac{\kappa_1}{s_n^2} \bar{u}_{ny} z_y dy dt 
+ \!\int_0^{T} \!\!
 (\frac{s_n'}{s_n} \bar{u}_n(\cdot,1) + \frac{1}{s_n} \psi(\bar{u}_n(\cdot,1))) z(\cdot,1) dt 
\nonumber \\
& = & \int_{Q(T)} (f_m(\bar{u}_n,\bar{v}_n) + \frac{s_n'}{s_n} y \bar{u}_{ny}) z dydt  
 \quad \mbox{ for } n.  \label{2-3q}
\end{eqnarray}
Since $s_n \to s$ in $C([0,T])$, from the above convergences it is clear that 
$$ \int_0^{T}\!\!\! \int_0^1 \bar{u}_{n_jt} z dx dt \to \int_0^T \langle \bar{u}_t, z \rangle_X dt,
\int_{Q(T)} \!\!\frac{\kappa_1}{s_{n_j}^2} \bar{u}_{n_jy} z_y dy dt \to 
\int_{Q(T)} \!\!\frac{\kappa_1}{s^2} \bar{u}_{y} z_y dy dt, 
$$
$$
\int_{Q(T)} (f_m(\bar{u}_{n_j},\bar{v}_{n_j}) + \frac{s_{n_j}'}{s_{n_j}} y \bar{u}_{n_jy}) z dydt 
\to  \int_{Q(T)} (f_m(\bar{u},\bar{v}) + \frac{s'}{s} y \bar{u}_{y}) z dydt
\mbox{ as } j \to \infty. $$
We show convergences of the third and fourth terms in the left hand side of (\ref{2-3q}) in the following way: 
\begin{eqnarray*}
& & | \int_0^{T}  (\frac{s_{n_j}'}{s_{n_j}} \bar{u}_{n_j}(\cdot,1) - \frac{s'}{s} \bar{u}(\cdot,1)) z(\cdot,1) dt |  \\
& \leq &  \int_0^{T}  |\frac{s_{n_j}'}{s_{n_j}}  - \frac{s'}{s}| 
  |\bar{u}_{n_j}(\cdot,1)|  |z(\cdot,1)| dt
+   \int_0^{T}  \frac{s'}{s} |\bar{u}_{n_j}(\cdot,1) -  \bar{u}(\cdot,1))| 
            |z(\cdot,1)| dt \\
& \leq &  u_* \int_0^{T}  |\frac{s_{n_j}'}{s_{n_j}}  - \frac{s'}{s}| |z|_X dt
+   \frac{\sqrt{2}}{s(0)}  \int_0^{T}  |s'| |\bar{u}_{n_j} -  \bar{u}|_H^{1/2} 
|\bar{u}_{n_jy} -  \bar{u}_y|_H^{1/2} |z|_H^{1/2} |z_y|_H^{1/2} dt, 
\end{eqnarray*}
and 
\begin{eqnarray*}
& & |\int_0^{T}  (\frac{1}{s_{n_j}} \psi(\bar{u}_{n_j}(\cdot,1)) - \frac{1}{s} \psi(\bar{u}(\cdot,1))) 
     z(\cdot,1) dt| \\
& \leq & \int_0^{T}  |\frac{1}{s_{n_j}} - \frac{1}{s}| 
          |\psi(\bar{u}_{n_j}(\cdot,1))| |z(\cdot,1)| dt 
+ \int_0^{T}  \frac{1}{s} |\psi(\bar{u}_{n_j}(\cdot,1)) - \psi(\bar{u}(\cdot,1))| 
     |z(\cdot,1)| dt \\
& \leq & \psi(u_*) \int_0^{T}  |\frac{1}{s_{n_j}} - \frac{1}{s}| |z|_X dt
+ \frac{C_5}{s(0)} \int_0^{T}   |\bar{u}_{n_j}(\cdot,1)) - (\bar{u}(\cdot,1))| |z(\cdot,1)| dt  \\
& \leq & \psi(u_*) \int_0^{T}  |\frac{1}{s_{n_j}} - \frac{1}{s}| |z|_X dt
+ \frac{\sqrt{2} C_5}{s(0)} \int_0^{T}   |\bar{u}_{n_j} - \bar{u}|_H^{1/2} 
      |\bar{u}_{n_jy} - \bar{u}_y|_H^{1/2} |z|_X dt  
     \mbox{ for } j, 
\end{eqnarray*}
where $C_5$ is a positive constant satisfying $|\psi(r) - \psi(r')| \leq C_5|r - r'|$ for $0 \leq r, r' \leq u_*$. Hence, we conclude that (SS3) holds. Note that we can get (SS4) in a similar fashion.

As next step, we prove the uniqueness of a weak solution to SP$_m(s, \bar{u}_0, \bar{v}_0$) on $[0,T]$. 
Let $(\bar{u}_1, \bar{v}_1)$ and $(\bar{u}_2, \bar{v}_2)$ be weak solutions of  SP$_m(s, \bar{u}_0, \bar{v}_0$) on $[0,T]$ and put $\bar{u} = \bar{u}_1 - \bar{u}_2$ and $\bar{v} = \bar{v}_1 - \bar{v}_2$ on $Q(T)$. Then (SS3) implies that 
\begin{eqnarray}
& & \langle \bar{u}_t, z \rangle_X + \frac{\kappa_1}{s^2} \int_0^1 \bar{u}_y z_y dy
         +  \frac{s'}{s} \bar{u}(\cdot,1) z(\cdot,1) 
    + \frac{1}{s} (\psi(\bar{u}_1(\cdot,1)) - \psi(\bar{u}_2(\cdot,1))  ) z(\cdot,1) 
   \nonumber \\ 
& = & \!\int_0^1  (f_m(\bar{u}_1,\bar{v}_1) - f_m(\bar{u}_2,\bar{v}_2) ) z dy 
+ \frac{s'}{s} \int_0^1 y \bar{u}_{y} z dy  
  \mbox{ for } z \in X \mbox{ a.e. on } [0,T].  \label{2-q3q}
\end{eqnarray}
By taking $z = \bar{u}$ in (\ref{2-q3q}) we have 
\begin{eqnarray*}
& & \frac{1}{2} \frac{d}{dt} 
         |\bar{u}|_H^2 + \frac{\kappa_1}{s^2} |\bar{u}_y|_H^2 
  + \frac{s'}{s} |\bar{u}(\cdot,1)|^2
 + \frac{1}{s} (\psi(\bar{u}_1(\cdot,1)) - \psi(\bar{u}_2(\cdot,1))) \bar{u}(\cdot,1) \\
& \leq &  C_m (|\bar{u}|_H + |\bar{v}|_H) |\bar{u}|_H + \frac{\kappa_1}{2s^2} |\bar{u}_y|_H^2 
+ \frac{1}{2\kappa_1} |s'|^2 |\bar{u}|_H^2   \mbox{ a.e. on } [0,T]
\end{eqnarray*}
so that 
$$
\frac{1}{2}  \frac{d}{dt} 
       |\bar{u}|_H^2 + \frac{\kappa_1}{2s^2} |\bar{u}_y|_H^2 
 \leq   C_m (|\bar{u}|_H + |\bar{v}|_H) |\bar{u}|_H 
        + \frac{1}{2\kappa_1} |s'|^2 |\bar{u}|_H^2   \mbox{ a.e. on } [0,T]. 
$$
We can also obtain the  inequality for $\bar{v}$. Accordingly, 
by adding these two inequalities and 
Gronwall's inequality we show the uniqueness. 

Finally, in order to prove (\ref{2x}), we put $u(t,x) = \bar{u}(t,\frac{x}{s(t)})$ and  $v(t,x) = \bar{v}(t,\frac{x}{s(t)})$ for $(t,x) \in Q_s(T)$ and 
$u_n(t,x) = \bar{u}_n(t,\frac{x}{s_n(t)})$ and  $v_n(t,x) = \bar{v}_n(t,\frac{x}{s_n(t)})$ for $(t,x) \in Q_{s_n}(T)$ and $n$. Then Lemma \ref{ineq} guarantees the following inequality: 
\begin{eqnarray}
& & \frac{1}{2} \frac{d}{dt} \int_0^{s_n} |u_n - g|^2 dx 
+  \frac{\gamma}{2} \frac{d}{dt} \int_0^{s_n} |v_n - h|^2 dx \nonumber \\
& & + \kappa_1 \int_0^{s_n} |u_{nx}|^2 dx 
  + \kappa_2 \gamma \int_0^{s_n} |v_{nx}|^2 dx 
+ \psi(u_n(\cdot,s_n)) (u(\cdot,s_n) - g)    \nonumber
 \\
& & + \frac{1}{2} s_n'( |u_n(\cdot,s_n) - g|^2 +  \gamma |v_n(\cdot,s_n) - h|^2) 
+ C_{\phi} \int_0^{s_n} |\gamma v_n - u_n|^{q + 1} dx \nonumber \\
& \leq & 
- g' \int_0^{s_n} (u_n - g)dx 
- \gamma h' \int_0^{s_n} (v_n - h)dx \label{2kk} \\
& & -  \int_0^{s_n} \phi_m (\gamma v_n - u_n)  (g - g_*) dx  
+ \gamma \int_0^{s_n} \phi_m (\gamma v_n - u_n)  (h - h_*) dx \nonumber \\
& & - s_n' g (u_n(\cdot,s_n) - g) - \gamma s_n' h (v_n(\cdot,s_n) - h) 
\quad \mbox{  a.e. on } [0,T], \nonumber
\end{eqnarray}
We integrate (\ref{2kk}) on $[0,t_1]$ with respect to $t$ for $0 < t_1 \leq T$. Then on account of the lower semi continuity of integral, we obtain  by letting $n \to  \infty$ 
\begin{eqnarray}
& & \frac{1}{2}  \int_0^{s(t_1)} |u(t_1) - g(t_1)|^2 dx 
+  \frac{\gamma}{2}  \int_0^{s(t_1)} |v(t_1) - h(t_1)|^2 dx \nonumber \\
& & + \kappa_1 \int_0^{t_1} \int_0^{s} |u_{x}|^2 dx dt 
  + \kappa_2 \gamma \int_0^{t_1} \int_0^{s} |v_{x}|^2 dx dt  
+  \int_0^{t_1} \psi(u(\cdot,s))(u(\cdot,s) - g) dt    \nonumber
 \\
& & + \frac{1}{2}  \int_0^{t_1} s'( |u(\cdot,s) - g|^2 +  \gamma |v(\cdot,s) - h|^2) dt  
+ C_{\phi}  \int_0^{t_1} \int_0^{s} |\gamma v - u|^{q + 1} dx dt \nonumber \\
& \leq &  
-  \int_0^{t_1} g' \int_0^{s} (u - g) dx dt 
-   \int_0^{t_1} \gamma h' \int_0^{s} (v - h)dx dt \nonumber \\
& & -   \int_0^{t_1} \int_0^{s} \phi_m (\gamma v - u)  (g - g_*) dx dt 
+ \gamma   \int_0^{t_1} \int_0^{s} \phi_m (\gamma v - u)  (h - h_*) dx dt \nonumber \\
& & -  \int_0^{t_1} (s' g (u(\cdot,s) - g) + \gamma s' h (v(\cdot,s) - h)) dt
\quad \mbox{ for }  0 < t_1 \leq T.  \nonumber
\end{eqnarray}
Relying on uniqueness, $(\bar{u}, \bar{v})$ is also a weak solution of the problem SP$_m(s, \bar{u}(t_0), \bar{v}(t_0)$) on $[t_0,T]$ for $0 < t_0 \leq T$. Hence, it holds that 
\begin{eqnarray}
& & \frac{1}{2}  \int_0^{s(t_1)} |u(t_1) - g(t_1)|^2 dx 
+  \frac{\gamma}{2}  \int_0^{s(t_1)} |v(t_1) - h(t_1)|^2 dx \nonumber \\
& & + \kappa_1 \int_{t_0}^{t_1} \int_0^{s} |u_{x}|^2 dx dt 
  + \kappa_2 \gamma \int_{t_0}^{t_1} \int_0^{s} |v_{x}|^2 dx dt  
+  \int_{t_0}^{t_1} \psi(u(\cdot,s))( u(\cdot,s) - g) dt    \nonumber
 \\
& & + \frac{1}{2}  \int_{t_0}^{t_1} s'( |u(\cdot,s) - g|^2 +  \gamma |v(\cdot,s) - h|^2) dt  
+ C_{\phi}  \int_{t_0}^{t_1} \int_0^{s} |\gamma v - u|^{q + 1} dx dt \nonumber \\
& \leq &  
-  \int_{t_0}^{t_1} g' \int_0^{s} (u - g) dx dt 
-   \int_{t_0}^{t_1} \gamma h' \int_0^{s} (v - h)dx dt \label{2xx} \\
& & -   \int_{t_0}^{t_1} \int_0^{s} \phi_m (\gamma v - u)  (g - g_*) dx dt 
+ \gamma   \int_{t_0}^{t_1} \int_0^{s} \phi_m (\gamma v - u)  (h - h_*) dx dt \nonumber \\ 
& & -  \int_{t_0}^{t_1} (s' g (u(\cdot,s) - g) + \gamma s' h(v(\cdot,s) - h))dt
\quad \mbox{ for }  0 \leq t_0 < t_1 \leq T.  \nonumber
\end{eqnarray}
Therefore, by dividing it by $t_1 - t_0$ and letting $t_1 \downarrow t_0$ we can obtain 
(\ref{2xx}). Thus we have proved this Proposition. \hfill $\Box$ 
 
\section{Interfaces propagate asymptotically like $\sqrt{t}$ as $t\to\infty$} 

In this section, we finally prove the main result -- Theorem \ref{main}.

\subsection{ Proof of the existence of a weak solution} 

We suppose (A1), (A2) and (A3). 
Then, since $f_m$ is Lipschitz continuous on $\R$ for each $m > 0$, 
by Theorem 1.1 of \cite{Ai-Mun1} 
P$(f_m)$ has a unique weak solution $\{s, u, v\}$ on $[0,T_m]$ for some $T_m > 0$. 

First, we show that P$(f_m)$ has a weak solution on $[0,\infty)$. In fact, let $[0,T_m^*)$ be the maximal interval of existence of a weak solution of  P$(f_m)$. We assume that $T_m^*$ is finite. 
Obviously, Lemma \ref{bounds} implies that 
$$ 0 \leq u \leq u_* \mbox{ and } 0 \leq v \leq v_* \mbox{ on } Q_s(T_m^*) $$
so that 
$s'(t) = \psi(u(t,s(t)) \leq \psi(u_*)$ for $0 \leq t < T_m^*$,
where $u_*$ and $v_*$ are positive constants given in the proof of Lemma \ref{bounds}. Accordingly, there exists a number $s(T_m^*) > 0$ such that $s(t) \to s(T_m^*)$ as $t \uparrow T_m^*$. Therefore, on account of (\ref{2x}) we infer that $\bar{u}, \bar{v} \in L^2(0,T_m^*; H^1(0,1))$, where 
$\bar{u}$ and $\bar{v}$ are functions defined by (\ref{trans}).
Similarly to (\ref{22}), 
$\bar{u}, \bar{v} \in W^{1,2}(0,T_m^*; X^*)$. This shows that there exist 
$\bar{u}(T_m^*), \bar{v}(T_m^*) \in L^{\infty}(0,1)$ such that 
$\bar{u}(t) \to \bar{u}(T_m^*)$ and $\bar{v}(t) \to \bar{v}(T_m^*)$ weakly in $H$ as 
$t \uparrow T_m^*$. Hence, by applying Theorem 1.1 of \cite{Ai-Mun1}, again, we can extend the solution beyond $T_m^*$. This is a contradiction, that is,  P$(f_m)$ has a weak solution on $[0,\infty)$. Moreover, it is obvious the weak solution of P$(f_m)$  is also a weak solution to 
P$(f)$, in case $m \geq \gamma v_* + u_*$. Thus we have proved the existence of a weak solution to P$(f)$ on $[0,\infty)$.  \hfill $\Box$

\subsection{ Proof of the upper estimate for the free boundary position} 
Let $\{s, u, v\}$ be a weak solution of P$(f)$ on $[0,\infty)$ and $\bar{u}$ and $\bar{v}$ are functions defined by (\ref{trans}). Then (S4) leads to:
\begin{eqnarray*}
&  &  \langle \bar{u}_t(t), z \rangle_X 
+  \frac{\kappa_1}{s^2(t)} \int_0^1 \bar{u}_y(t) z_y dy 
+  (\frac{s'(t)}{s(t)} \bar{u}(t,1) + \frac{s'(t)}{s(t)} )z(1) dt  \\
& = & \int_0^1 (f(\bar{u}(t),\bar{v}(t)) + \frac{s'(t)}{s(t)} y \bar{u}_y(t)) z
 dy \quad \mbox{ for } z \in X \mbox{ and a.e. } t \in [0,\infty).
\end{eqnarray*}

 Accordingly, by taking $z = s^2 y$,  we have 
\begin{eqnarray*}
& & \langle \bar{u}_t, s^2y \rangle_X + \kappa_1 \int_0^1 \bar{u}_y  dy
         +  s s' (\bar{u}(\cdot,1) + 1) 
   \nonumber \\ 
& = & \!\int_0^1  (f(\bar{u},\bar{v}) s^2y dy 
+ s s'  \int_0^1  \bar{u}_{y} y^2 dy  
  \quad \mbox{ a.e. on } [0,\infty).  \label{3-2}
\end{eqnarray*}
It is clear that (see \cite[Proposition 23.23]{Zeidler}) 
$$  \langle \bar{u}_t, s^2y \rangle_X = \frac{d}{dt}\int_0^1 \bar{u} s^2y dy 
                - \int_0^1 2 \bar{u} ss' y dy, \int_0^1 \bar{u} s^2y dy = \int_0^s x u dx  \mbox{ a.e. on } [0,\infty).  $$
It follows that 
\begin{equation}
  \frac{d}{dt}\int_0^s x u dx + \kappa_1 \int_0^1 \bar{u}_y  dy
             + s s'
  =  \!\int_0^1  f(\bar{u},\bar{v}) s^2y dy \mbox{ a.e. on } [0,\infty).  \label{3-3}\end{equation}
We can obtain the similar equation for $\bar{v}$ to (\ref{3-3}). Accordingly, we see that 
\begin{eqnarray}
 \frac{d}{dt}\int_0^s x (u  + v)  dx + \kappa_1 \int_0^1 \bar{u}_y  dy 
+ \kappa_2 \int_0^1 \bar{v}_y  dy + s s'   
 =  0 
   \mbox{ a.e. on } [0,\infty).  \label{3-19}
\end{eqnarray}
By integrating it, it holds that 
\begin{eqnarray*}
& &  \int_0^{s(t)} x (u(t)  + v(t))  dx + \kappa_1 \int_0^t u(\tau, s(\tau))  d\tau 
+ \kappa_2 \int_0^t v(\tau, s(\tau))  d\tau 
             + \frac{1}{2} s^2(t)  \\
&  = &  \int_0^{s_0} x (u_0  + v_0)  dx + \kappa_1 \int_0^t g(\tau)  d\tau 
+ \kappa_2 \int_0^t h(\tau)  d\tau 
             + \frac{1}{2} s^2_0 
 \quad   \mbox{  for } t \in [0,\infty).  
\end{eqnarray*}
Making use of the positivity of $u$ and $v$, we observe that 
$$ \frac{1}{2} s(t)^2 \leq \frac{1}{2} s_0^2 + \int_0^{s_0} x (u_0 + v_0) dx 
  + (\kappa_1 g^* + \kappa_2 h^*) t   \quad   \mbox{  for } t \in [0,\infty), 
$$
where $g^* = |g|_{L^{\infty}(0,\infty)}$ and 
$h^* = |h|_{L^{\infty}(0,\infty)}$. 
This inequality guarantees the existence of a  positive constant  $C_*$ satisfying
$$ s(t) \leq C_* \sqrt{t + 1} \mbox{ for } t \geq 0. \eqno{\Box} $$

\vskip 12pt
{\it Proof of the lower estimate for the free boundary. } 
First, we show 
\begin{equation}
 \int_0^t \int_0^s |v_x|^2 dx d\tau \leq K_1 (s(t)  + 1) \quad \mbox{ for } 
    t \geq 0,  \label{q67}
\end{equation}
where $K_1$ is a positive constant. In fact, Proposition \ref{prop2} implies 
\begin{eqnarray}
& & \frac{1}{2}  \int_0^{s(t)} |u(t) - g(t)|^2 dx 
+  \frac{\gamma}{2}  \int_0^{s(t)} |v(t) - h(t)|^2 dx \nonumber \\
& & + \kappa_1 \int_0^t \int_0^{s} |u_x|^2 dx d\tau 
  + \kappa_2 \gamma \int_0^t \int_0^{s} |v_x|^2 dx d\tau  
  +  \int_0^t s' u(\cdot,s) d\tau  
\nonumber
 \\
& \leq & 
\frac{1}{2} \int_0^{s_0} |u_0 - g(0)|^2 dx +
\frac{\gamma}{2} \int_0^{s_0} |v_0 - h(0)|^2 dx \nonumber \\
& & + \int_0^t s' g d\tau  
- \int_0^t g' \int_0^{s} (u - g)dx d\tau 
- \gamma  \int_0^t h' \int_0^{s} (v - h)dx d\tau  \nonumber  \\
& & -  \int_0^t \int_0^{s} \phi (\gamma v - u)  (g - g_*) dx d\tau  
+ \gamma \int_0^t \int_0^{s} \phi (\gamma v - u)  (h - h_*) dx d\tau  \nonumber \\
& & - \int_0^t (s' g (u(\cdot,s) - g) + \gamma s' h (v(\cdot,s) - h)) d\tau 
 \nonumber \\
& \leq & 
\frac{1}{2} \int_0^{s_0} |u_0 - g(0)|^2 dx +
\frac{\gamma}{2} \int_0^{s_0} |v_0 - h(0)|^2 dx \label{3q} \\
& & + g^*( s(t) - s_0)   
+ (u_* + g^*) s(t) \int_0^t |g'|  d\tau 
+ \gamma (v_* + h^*) s(t)  \int_0^t |h'| d\tau  \nonumber  \\
& & + s(t)  C_4  (\int_0^t  |g - g_*| dx d\tau  
+ \gamma \int_0^t |h - h_*| dx d\tau)  \nonumber \\
& & + (g^*(u_** + g^*)  + \gamma  h^* (v_* + h^*)) (s(t) - s_0)
\quad \mbox{  for }  t \geq 0. \nonumber
\end{eqnarray}
Obviously, by (A2) we can take a positive number $K_1$ satisfying (\ref{q67}). 

Recalling (\ref{3-19}), we have 
\begin{eqnarray*}
& &  \int_0^{s(t)} x (u(t)  + v(t))  dx + \kappa_1 \int_0^t u(\tau, s(\tau))  d\tau 
+ \kappa_2 \int_{Q_s(t)} v_x  dx d\tau 
             + \frac{1}{2} s^2(t)  \\
&  = &  \int_0^{s_0} x (u_0  + v_0)  dx + \kappa_1 \int_0^t g(\tau)  d\tau 
             + \frac{1}{2} s^2_0  \\
&  \geq  &   \kappa_1 g_0 t
 \quad   \mbox{  for } t\geq 0  
\end{eqnarray*}
so that 
\begin{eqnarray*}
   \kappa_1 g_0 t 
& \leq & \kappa_2 (\int_{Q_s(t)} |v_x|^2 dx d\tau)^{1/2} (s(t)t)^{1/2} 
  + (u_* + v_*) \int_0^{s(t)} x dx  +   \frac{1}{2} s^2(t) \\ 
& &  + \frac{\kappa_1}{\kappa_0^{1/p}} \int_0^t |s'|^{1/p} d\tau \\
& \leq & \kappa_2 (K_1 s(t) + K_1)^{1/2} (s(t)t)^{1/2} 
  + \frac{1}{2} (u_* + v_*) s(t)^2  +   \frac{1}{2} s^2(t) \\ 
& &  + \frac{\kappa_1}{\kappa_0^{1/p}} (\int_0^t |s'| d\tau)^{1/p} t^{1-1/p} \\
& \leq & K_2 (s(t) + 1)s(t) + \frac{\kappa_1 g_0}{2} t 
 \quad   \mbox{  for } t \geq 0,  
\end{eqnarray*}
where $K_2$ is a positive constant. 
Then it holds that 
$$ \frac{\kappa_1 g_0}{2} t \leq K_2 s(t)^2 + \frac{\kappa_1 g_0}{4} 
  + \frac{1}{\kappa_1 g_0} K_2^2 s(t)^2  \mbox{  for } t \in [0,T]. $$
Hence, it is easy to see that 
$$ \frac{\kappa_1 g_0}{4} t \leq (K_2   
  + \frac{1}{\kappa_1 g_0} K_2^2) s(t)^2  \mbox{  for } t \geq 1. $$
In case $0 \leq t \leq 1$, we have 
$$ s_0 \sqrt{t} \leq s_0 \leq s(t). $$
Thus we have now completed the proof of the Theorem.  \hfill $\Box$

%\bibliographystyle{elsarticle-num}
%\bibliographystyle{plain}
%\bibliography{FbpLarge}

\begin{thebibliography}{10}
\expandafter\ifx\csname url\endcsname\relax
  \def\url#1{\texttt{#1}}\fi
\expandafter\ifx\csname urlprefix\endcsname\relax\def\urlprefix{URL }\fi
\expandafter\ifx\csname href\endcsname\relax
  \def\href#1#2{#2} \def\path#1{#1}\fi

\bibitem{Has}
L.~Haselbach, Potential for carbon dioxide absorption in concrete, Journal of
  Environmental Engineering 135 (2009) 465--472.

\bibitem{Ruan}
X.~Ruan, Z.~Pan, Mesoscopic simulation method of concrete carbonation process,
  Structure and Infrastructure Engineering 8~(2) (2012) 99--110.

\bibitem{Ai-Mun1}
T.~Aiki, A.~Muntean, Existence and uniqueness of solutions to a mathematical
  model predicting service life of concrete structures, Adv. Math. Sci. Appl.
  19 (2009) 109--129.

\bibitem{Malte}
M.~A. Peter, M.~B{\"o}hm, Different choices of scaling in homogenization of
  diffusion and interfacial exchange in a porous medium, Math. Meth. Appl. Sci.
  31~(11) (2008) 1257 -- 1282.

\bibitem{Maria}
A.~Muntean, M.~Neuss-Radu, A multiscale {G}alerkin approach for a class of
  nonlinear coupled reaction-diffusion systems in complex media, Journal of
  Mathematical Analysis and Applications 371~(2) (2010) 705 -- 718.

\bibitem{IFB}
T.~Aiki, A.~Muntean, A free-boundary problem for concrete carbonation: rigorous
  justification of the $\sqrt t$-law of propagation, Interfaces and Free
  Boundaries 15 (2013) in press.

\bibitem{CPAA}
T.~Aiki, A.~Muntean, Large time behavior of solutions to a moving-interface
  problem modeling concrete carbonation, Communications on Pure and Applied
  Analysis 9 (2010) 1117--1129.

\bibitem{Kris}
K.~Sisomphon, Influence of {P}ozzolanic {M}aterial {A}dditions on the
  {D}evelopment of the {A}lkalinity and the {C}arbonation {B}ehavior of
  {C}omposite {C}ement {P}astes and {C}oncretes, Ph.D. thesis, TU
  Hamburg-Harburg, Germany (2004).

\bibitem{Lions}
J.~L. Lions, Quelques m\'ethodes de resolution des probl\`emes aux limites
  non-lin\'eaires, Dunod, Paris, 1990.

\bibitem{Zeidler}
E.~Zeidler, Nonlinear {F}unctional {A}nalysis and its {A}pplications. {L}inear
  {M}onotone {O}perators, II/A, Springer Verlag, NY, Berlin, 1969.

\end{thebibliography}

%\medskip
%Received xxxx 20xx; revised xxxx 20xx.
%\medskip

\end{document}